# ESTIMATION IN SPIN GLASSES: A FIRST STEP[1]

By Sourav Chatterjee

*University of California, Berkeley*

The Sherrington–Kirkpatrick model of spin glasses, the Hopfield model of neural networks and the Ising spin glass are all models of binary data belonging to the one-parameter exponential family with quadratic sufficient statistic. Under bare minimal conditions, we establish the $\sqrt{N}$-consistency of the maximum pseudolikelihood estimate of the natural parameter in this family, even at critical temperatures. Since very little is known about the low and critical temperature regimes of these extremely difficult models, the proof requires several new ideas. The author's version of Stein's method is a particularly useful tool. We aim to introduce these techniques into the realm of mathematical statistics through an example and present some open questions.

**1. Introduction and main result.**

1.1. *Statement of the problem.* Suppose our data is a vector of dependent $\pm 1$-valued random variables, denoted by $\sigma = (\sigma_1, \ldots, \sigma_N)$. A simple and natural way of modeling the dependence is to consider the one-parameter exponential family where the sufficient statistic is a quadratic form. Explicitly, this means that we have a collection of real numbers $(J_{ij}^N)_{1 \le i < j \le N}$ defining a parametric family $(\mathbb{P}_\beta^N)_{\beta \ge 0}$ of probability distributions on $\{-1, 1\}^N$ in the following way: For any $\tau \in \{-1, 1\}^N$ and any $\beta \ge 0$,

$$\mathbb{P}_\beta^N\{\sigma = \tau\} = 2^{-N} \exp(\beta H_N(\tau) - N\psi_N(\beta)), \tag{1}$$

where

$$H_N(\tau) := \sum_{1 \le i < j \le N} J_{ij}^N \tau_i \tau_j. \tag{2}$$

Received April 2006; revised November 2006.
[1]Supported in part by NSF Grant DMS-07-07054 and a Sloan Research Fellowship.
*AMS 2000 subject classifications.* 62F10, 62F12, 60K35, 82B44.
*Key words and phrases.* Spin glass, neural networks, estimation, consistency, exponential families.

Here $\beta$ is the only unknown parameter. The function $\psi_N$ is determined by the normalizing condition $\sum_\tau \mathbb{P}^N_\beta\{\sigma = \tau\} = 1$. For future use, we let $J^N_{ji} = J^N_{ij}$ and $J^N_{ii} = 0$ for each $i$ and $j$ and let $J^N$ be the matrix $(J^N_{ij})_{1 \leq i,j \leq N}$.

Many popular models, for example, the usual ferromagnetic Ising model [18, 21], the Sherrington–Kirkpatrick mean field model [24] of spin glasses and the Hopfield model [17] of neural networks all belong to the above family, with various special structures for $J^N$.

The main problem with the maximum likelihood estimation of $\beta$ is that the function $\psi_N$ is generally not computable, theoretically or otherwise. There are numerical methods for computing approximate MLEs (e.g., Geyer and Thompson [12]), but very little is known about the number of steps required for convergence. Moreover, we do not know of any gradient-based algorithm which provably converges to a *global* maximum of the likelihood function, and there are serious doubts whether that is indeed possible. The Jerrum–Sinclair [20] algorithm for computing the normalizing constant converges in polynomial time, but requires the presence of a local dependency graph, which makes it inapplicable to spin glasses. Most importantly, even if one assumes that the MLE has somehow been approximated, general conditions for the consistency of the MLE are not available.

1.2. *The pseudolikelihood estimator.* A natural recourse is to consider Julian Besag's *maximum pseudolikelihood estimator* (MPLE) [4, 5]. In short, it gives the following prescription: Suppose we have a random vector $X = (X_1, \ldots, X_n)$ whose distribution is parametrized by a parameter $\beta$. Let $f_i(\beta, X)$ be the conditional probability density of $X_i$ given $(X_j)_{j \neq i}$. Then the MPLE of $\beta$ is defined as

$$\hat{\beta}_{MPL} := \arg\max_\beta \prod_i f_i(\beta, X).$$

In our setting, the conditional density of $\sigma_i$ given the rest is easy to compute:

$$\log f_i(\beta, \sigma) = \beta \sigma_i \sum_{j=1}^N J^N_{ij} \sigma_j - \log \cosh\left(\beta \sum_{j=1}^N J^N_{ij} \sigma_j\right) - \log 2.$$

Thus, the derivative of the pseudolikelihood function is

$$\begin{aligned}
S_\sigma(\beta) &:= \frac{\partial}{\partial \beta} \sum_{i=1}^N \log f_i(\beta, \sigma) \\
&= \sum_{i,j=1}^N J^N_{ij} \sigma_i \sigma_j - \sum_{i,j=1}^N J^N_{ij} \sigma_j \tanh\left(\beta \sum_{k=1}^n J^N_{ik} \sigma_k\right).
\end{aligned} \tag{3}$$

Noting that $H_N(\sigma) = \frac{1}{2}\sum_{i,j=1}^{N} J_{ij}^N \sigma_i \sigma_j$, we use the above expression to define the maximum pseudolikelihood estimate in this problem as

$$(4) \qquad \hat{\beta}_N := \inf\left\{x \geq 0 : H_N(\sigma) = \tfrac{1}{2}\sum_{i,j=1}^{N} J_{ij}^N \sigma_j \tanh\left(x \sum_{k=1}^{N} J_{ik}^N \sigma_k\right)\right\}.$$

The infimum of the empty set is defined to be $\infty$, as usual. It is not difficult to show that the expression on the right-hand side is an increasing function of $x$. Thus, it is not only feasible but extremely easy to compute $\hat{\beta}_N$, either by Newton–Raphson or even a simple grid search.

1.3. *The consistency result.* Having defined the MPLE, it is natural to ask whether it is consistent. Unlike other exponential families, spin glass models present the unique challenge that even the most basic characteristics, like the correlations between spins at different sites—let alone weak laws or central limit theorems—are completely intractable when $\beta$ is larger than a threshold. Any attempt at proving a consistency result in these so-called "low temperature regimes" has to overcome the lack of almost any meaningful information about the behavior of the model. For instance, in most models of spin glasses, there is no known technique for proving that the function $S_\sigma$ defined in (3) converges to a nonrandom limit function (after suitable normalization) as $N \to \infty$, unless $\beta$ is sufficiently small. The main achievement of this paper is a new technique that bypasses these hurdles and proves, with minimal information, that the MPLE is $\sqrt{N}$-consistent at all temperatures. The most challenging part is to tackle the issues at the "critical temperatures," that is, points of phase transition.

Recall the following definition from basic linear algebra: The $L^2$ operator norm of a square matrix $A$ is defined as $\|A\| := \sup_{\|x\|=1}\|Ax\|$, where $\|x\|$ stands for the usual Euclidean norm of the vector $x$. If $A$ is a real symmetric matrix, then $\|A\|$ is equal to the spectral radius of $A$. Our main result is the following.

THEOREM 1.1. *Consider the exponential family of models (1) with the quadratic sufficient statistic (2). Fix $\beta > 0$, and let $\hat{\beta}_N$ be the estimator (4). Suppose we have a sequence of such models with $N \to \infty$, satisfying*

(a) $\sup_N \|J^N\| < \infty$, *and*
(b) $\liminf_{N\to\infty} \psi_N(\beta) > 0$.

*Then $\{\hat{\beta}_N\}_{N \geq 1}$ is a $\sqrt{N}$-consistent sequence of estimators for $\beta$.*

Note that $\psi_N(\beta)$ is always an increasing nonnegative function, as long as $\beta > 0$. Hence, if condition (b) is satisfied for some positive $\beta$, then it holds for every $\beta' \geq \beta$.

Surprisingly, the MPLE may not be consistent at $\beta = 0$. A natural counterexample, which also shows that $\hat{\beta}_N$ may not be consistent if condition (b) is not satisfied, is provided by the Curie–Weiss model of ferromagnetic interaction at high temperature. This is discussed in Section 1.7 at the end of this section.

Some rigorous results are known about the consistency of the MPLE in settings with some kind of Markovian structure (e.g., lattice processes [10, 13] and spatial point processes [19]), but the techniques of these papers cannot be used to prove Theorem 1.1. The reason is that the dependence in spin glasses is neither local nor mean field in the classical sense, and there is no way to extract any conditional independence.

The main technique employed in this paper is a new version of Charles Stein's method of exchangeable pairs, developed by the author in [7, 8, 9]. The first step is the construction of an exchangeable pair $(\sigma, \sigma')$ and an antisymmetric function $F$ such that $\mathbb{E}(F(\sigma, \sigma')|\sigma) =$ the pseudolikelihood score function. The theory then gives a "temperature-free" method of showing that the derivative (3) of the log-pseudolikelihood function is close to zero with high probability at the true value of the parameter, resulting in the following lemma.

LEMMA 1.2. *Let $S_\sigma$ be the derivative of the log-pseudolikelihood defined in (3). Then for any $\beta \geq 0$ we have*

$$\mathbb{E}_\beta(S_\sigma(\beta)^2) \leq \frac{C}{N},$$

*where $C = 6\|J^N\|^2 + 6\beta\|J^N\|^3 + 2\beta^2\|J^N\|^4$. Consequently, $\mathbb{P}_\beta\{|S_\sigma(\beta)| > \delta\} \leq C/(N\delta^2)$ for any $\delta > 0$.*

Having shown that, the next step is to prove consistency of the MPLE by inverting the pseudolikelihood function. This is usually an easy step in theoretical statistics; here, however, this is a very hard step because the most basic features of the model are unknown. To see that there are nontrivial issues involved, one can just look at the counterexample of the Curie–Weiss model at high temperature, presented in Section 1.7, where Lemma 1.2 holds but the MPLE is not consistent.

Theorem 1.1 has the desirable feature that the conditions (a) and (b) are satisfied in almost all commonly used models of spin glasses. We consider two examples.

1.4. *Application to the Sherrington–Kirkpatrick (S–K) model.* In the Sherrington–Kirkpatrick model of spin glasses [24], we have

$$J_{ij}^N = N^{-1/2} g_{ij},$$

where $(g_{ij})_{1 \leq i < j \leq \infty}$ is a fixed realization of an array of independent standard Gaussian random variables (and $g_{ji} = g_{ij}$). Some of the most important questions about the S–K model have only recently been answered, mainly due to the monumental efforts of Talagrand (see [25, 26], to mention a few) and the important contributions of Guerra [14], Guerra and Toninelli [15], Panchenko [22], Comets and Neveu [11] and other authors (e.g., [1, 2, 23]).

However, despite all the progress, the model is still far from tractable, especially when $\beta \geq 1$. Thus it is quite remarkable that Theorem 1.1 applies with equal ease for all $\beta$, as we now show. It follows from a standard result in random matrix theory (see, e.g., Bai [3], Theorem 2.12) that condition (a) holds for almost every realization of $(g_{ij})_{i,j}$. For (b), we can use the monotonicity of $\psi_N$, and the result [1] that $\lim \psi_N(\beta) = \beta^2/4$ for $\beta < 1$ in the S–K model (see Talagrand [25], Theorem 2.2.1), to conclude that $\liminf \psi_N(\beta) > 0$ for every positive $\beta$.

1.5. *Application to the Hopfield model.* In the Hopfield model [17] for a system with $N$ particles affected by $M$ attractors, we have

$$J_{ij}^{M,N} = \frac{1}{N} \sum_{k=1}^{M} \eta_{ik} \eta_{jk},$$

where $\{\eta_{ik}\}_{i \leq N, k \leq M}$ is a fixed realization of a collection of independent random variables with $\mathbb{P}\{\eta_{ik} = \pm 1\} = 1/2$. This model is even less understood than the S–K model. The groundbreaking contributions in the study of this model are mainly due to Bovier and Gayrard (see [6]) and Talagrand (compiled in [25], Chapter 5).

Again, it follows from random matrix theory (Bai [3], Section 2.2.2) that condition (a) is almost surely satisfied whenever $M/N$ stays bounded as $N \to \infty$. The validity of (b) follows from spin glass theory (Talagrand [25], Theorem 5.2.1). The conditions can similarly be verified for the other models mentioned in the abstract.

1.6. *Multiparameter models and open problems.* The main shortcoming of Theorem 1.1 is that it applies only to one-parameter families. The good news is that the main tool, Lemma 1.2, can easily extend to multiparameter families. However, as mentioned before, the surprising fact is that proving consistency of the $\hat{\beta}_N$ is very hard even after one has something like Lemma 1.2. The argument that we use in this paper to move from Lemma 1.2 to the consistency of $\hat{\beta}_N$ is intensely one-dimensional. It is so specialized that it breaks down even if we just add a simple linear term like $h \sum_{i=1}^{n} \sigma_i$ in the Hamiltonian.

The main open problems (currently under investigation by the author) are (i) to prove a version of Lemma 1.2 for multiparameter models, which should be doable by the techniques of this paper, and (ii) to deduce a multi-parameter version of Theorem 1.1 from this generalized form of Lemma 1.2, which is likely to be quite hard if not impossible; the author does not know how to solve this part with the current tools.

Besides these issues, there are also natural open questions about variance bounds and central limit theorems for the pseudolikelihood estimates, with the tantalizing possibility that non-Gaussian limits hold at the critical temperatures.

1.7. *A counterexample.* Theorem 1.1 may fail if condition (b) is not satisfied. In particular, this may happen at $\beta = 0$. Counterexamples are easy to construct if (a) is not satisfied; necessity of (b) is more subtle.

The simplest counterexample is provided by the Curie–Weiss model where $J_{ij}^N = 1/N$ for all $i, j$. This is a well-known toy model of ferromagnetic interaction. It is known that for $0 \le \beta < 1$, $\lim_{N \to \infty} \psi_N(\beta) = 0$ in this model (see, e.g., [25], page 324). Let $m_N = \frac{1}{N} \sum_{i=1}^{N} \sigma_i$, and define the function

$$f_N(x) := \frac{S_\sigma(x)}{2H_N(\sigma)} = 1 - \frac{\tanh(xm_N)}{m_N}.$$

Suppose $0 \le \beta < 1$. It is known that $m_N \to 0$ in probability in this regime (again, see [25], page 324). Thus, $f_N(x) \to 1 - x$ in probability for each $x$. Since $f_N$'s are increasing functions, and $\hat{\beta}_N$ is a zero of $f_N$, a standard argument now gives that $\hat{\beta}_N \to 1$ in probability. Thus, the pseudolikelihood estimate of $\beta$ in the Curie–Weiss model is not consistent when $0 \le \beta < 1$. Note that Lemma 1.2 continues to hold, though.

**2. A finite sample result and proofs.** In this section $N$ is fixed. Consequently, subscripts and superscripts involving $N$ are unnecessary. Let $J$ be an $N \times N$ nonzero symmetric matrix with zeros on the diagonal. As before, consider the one-parameter exponential family of probability distributions $(\mathbb{P}_\beta)_{\beta \ge 0}$ on $\{-1, 1\}^N$, defined as

$$\mathbb{P}_\beta(\{\tau\}) = 2^{-N} \exp(\beta H(\tau) - N\psi(\beta)),$$

where the sufficient statistic $H$ has the form

$$H(\tau) = \sum_{1 \le i < j \le N} J_{ij} \tau_i \tau_j.$$

Here $J$ is an $N \times N$ nonzero symmetric matrix with zeros on the diagonal. As usual $\tau$ will always denote a typical element of $\{-1, 1\}^N$ and $\sigma$ will be reserved for random elements. Now for each $i$, let us define the function $m_i : \{-1, 1\}^N \to \mathbb{R}$ as

(5) $$m_i(\tau) = \sum_{1 \le j \le N} J_{ij} \tau_j.$$

Note that $m_i(\tau)$ does not depend on $\tau_i$ because $J_{ii} = 0$ by assumption. It is not difficult to verify that if $\sigma \sim \mathbb{P}_\beta$, then

$$\mathbb{E}_\beta(\sigma_i|(\sigma_j)_{j\ne i}) = \tanh(\beta m_i(\sigma)).$$

The quantity $m_i(\tau)$ is called the *local field at site $i$* in the configuration $\tau$. Note that $H(\tau) = \frac{1}{2}\sum_i m_i(\tau)\tau_i$ in the class of models that we are considering. Next, for each $\tau \in \{-1,1\}^N$, we define a map $S_\tau : [0,\infty) \to \mathbb{R}$ as

$$(6) \qquad S_\tau(x) = \frac{1}{N}\sum_{1 \le i \le N} m_i(\tau)(\tau_i - \tanh(xm_i(\tau))).$$

Note that this is the same as the derivative of the log-pseudolikelihood defined in (3). Interpreting $\tanh(\pm\infty)$ as $\pm 1$, we can continuously extend $S_\tau$ to $[0,\infty]$ by defining

$$S_\tau(\infty) = \frac{1}{N}\sum_{i=1}^N (m_i(\tau)\tau_i - |m_i(\tau)|).$$

Note that $S_\tau(\infty) \le 0$ for all $\tau$. Finally, let

$$(7) \qquad \hat{\beta}(\tau) := \inf\{x \ge 0 : S_\tau(x) = 0\}.$$

It is easy to check that this is exactly the estimate defined in (4). When $\sigma$ is a random element picked from the measure $\mathbb{P}_\beta$, we will usually just write $\hat\beta$ instead of $\hat\beta(\sigma)$.

THEOREM 2.1. *Take any $\beta > 0$ and $0 < \varepsilon < 1$. Let $\gamma = \max\{\beta\|J\|, 1\}$. Let $\hat\beta$ be the estimate defined above in (7). Then*

$$\mathbb{P}_\beta\{|\tanh(C\hat\beta) - \tanh(C\beta)| > \varepsilon\} \le \frac{K}{N\varepsilon^2},$$

*where we can take*

$$C = \frac{8\beta\|J\|^2}{\psi(\beta)} \quad \text{and} \quad K = \frac{8\gamma^2}{\psi(\beta)^3} + \frac{2^{30}\gamma^{12}}{9\psi(\beta)^{10}}.$$

Note that Theorem 1.1 follows directly from Theorem 2.1: Suppose the conditions (a) and (b) of Theorem 1.1 are satisfied. By combining Theorem 2.1 with conditions (a) and (b), we get

$$\tanh(C_N\hat\beta_N) - \tanh(C_N\beta) \to 0 \qquad \text{in probability,}$$

where $C_N = 8\beta\|J^N\|^2/\psi_N(\beta)$. A simple verification shows that $\psi_N(\beta) \le \beta\|J^N\|$. Thus, by condition (b), we get $\liminf \|J^N\| > 0$. Combined with (a), this gives $0 < \liminf C_N \le \limsup C_N < \infty$. It follows that $\hat\beta_N - \beta = O(N^{-1/2})$.

Let us now begin the proof of Theorem 2.1 with the observation that the model remains unchanged if we replace $\beta$ by $\beta\|J\|$ and $J$ by $\|J\|^{-1}J$. Thus, without loss of generality, we can assume that

(8) $$\|J\| = 1.$$

The proof is divided into a sequence of lemmas. There are two main steps: (1) To show that $S_\sigma(\beta) \approx 0$ with high probability under $\mathbb{P}_\beta$ (this is Lemma 1.2), and (2) to show that $S_\sigma(\beta) \approx 0 \Rightarrow \beta \approx \hat\beta$. While the first step is the conceptual mainstay, the second step is surprisingly hard because of the extreme lack of information about the model. It is carried out in four steps (Lemmas 2.2, 2.3, 2.4 and 2.5).

PROOF OF LEMMA 1.2. Recall that we are working under the assumption that $\|J\| = 1$. Define a function $F: \{-1,1\}^N \times \{-1,1\}^N \to \mathbb{R}$ as

$$F(\tau, \tau') = \tfrac{1}{2} \sum_{i=1}^{N} (m_i(\tau) + m_i(\tau'))(\tau_i - \tau_i').$$

Note that $F$ is *antisymmetric*, that is, $F(\tau, \tau') \equiv -F(\tau', \tau)$. This will be useful later on.

Now fix $\beta \geq 0$ and suppose $\sigma$ is drawn from the Gibbs measure at inverse temperature $\beta$. Since $\beta$ is fixed in this proof, we will write $\mathbb{E}$ and $\mathbb{P}$ instead of $\mathbb{E}_\beta$ and $\mathbb{P}_\beta$.

Now choose a coordinate $I$ uniformly at random, and replace the $I$th coordinate of $\sigma$ by a sample drawn from the conditional distribution of $\sigma_I$ given $(\sigma_j)_{j \neq I}$. Call the resulting vector $\sigma'$. Then $(\sigma, \sigma')$ is an exchangeable pair of random variables. Observe that

$$F(\sigma, \sigma') = m_I(\sigma)(\sigma_I - \sigma_I').$$

Now let

$$f(\sigma) := \mathbb{E}(F(\sigma, \sigma')|\sigma) = \frac{1}{N} \sum_{i=1}^{N} m_i(\sigma)(\sigma_i - \mathbb{E}(\sigma_i|(\sigma_j)_{j \neq i}))$$

$$= \frac{1}{N} \sum_{i=1}^{N} m_i(\sigma)(\sigma_i - \tanh(\beta m_i(\sigma))) = S_\sigma(\beta).$$

Then $\mathbb{E}(f(\sigma)^2) = \mathbb{E}(f(\sigma)F(\sigma, \sigma'))$. Since $(\sigma, \sigma')$ is an exchangeable pair, therefore

$$\mathbb{E}(f(\sigma)F(\sigma, \sigma')) = \mathbb{E}(f(\sigma')F(\sigma', \sigma)).$$

Again, since $F$ is antisymmetric, we have $\mathbb{E}(f(\sigma')F(\sigma', \sigma)) = -\mathbb{E}(f(\sigma')F(\sigma, \sigma'))$. Combining, we have

(9)
$$\mathbb{E}(f(\sigma)^2) = \mathbb{E}(f(\sigma)F(\sigma, \sigma')) = -\mathbb{E}(f(\sigma')F(\sigma, \sigma'))$$
$$= \tfrac{1}{2}\mathbb{E}((f(\sigma) - f(\sigma'))F(\sigma, \sigma')).$$

For any $1 \leq j \leq N$ and $\tau \in \{-1,1\}^N$, let
$$\tau^{(j)} := (\tau_1, \ldots, \tau_{j-1}, -\tau_j, \tau_{j+1}, \ldots, \tau_N)$$
and
$$p_j(\tau) := \frac{e^{-\beta \tau_j m_j(\tau)}}{e^{\beta m_j(\tau)} + e^{-\beta m_j(\tau)}} = \mathbb{P}\{\sigma'_j = -\tau_j | \sigma = \tau, I = j\}.$$
Then
$$\mathbb{E}((f(\sigma) - f(\sigma'))F(\sigma, \sigma')|\sigma)$$
(10)
$$= \frac{1}{N} \sum_{j=1}^N (f(\sigma) - f(\sigma^{(j)})) F(\sigma, \sigma^{(j)}) p_j(\sigma)$$
$$= \frac{1}{N} \sum_{j=1}^N (f(\sigma) - f(\sigma^{(j)}))(2 m_j(\sigma)\sigma_j) p_j(\sigma).$$

Now, for ease of notation, we define the functions $a_i$ and $b_{ij}$ as
(11)
$$a_i(\tau) := \tau_i - \tanh(\beta m_i(\tau))$$
and
(12)
$$b_{ij}(\tau) := \tanh(\beta m_i(\tau)) - \tanh(\beta m_i(\tau^{(j)})).$$
Then $f(\tau) = N^{-1} \sum_i m_i(\tau) a_i(\tau)$, and hence
$$f(\sigma) - f(\sigma^{(j)})$$
$$= \frac{1}{N} \sum_{i=1}^N (m_i(\sigma) - m_i(\sigma^{(j)})) a_i(\sigma) + \frac{1}{N} \sum_{i=1}^N m_i(\sigma^{(j)})(a_i(\sigma) - a_i(\sigma^{(j)}))$$
$$= \frac{2\sigma_j}{N} \sum_{i=1}^N J_{ij} a_i(\sigma) + \frac{2 m_j(\sigma)\sigma_j}{N} - \frac{1}{N} \sum_{i=1}^N m_i(\sigma^{(j)}) b_{ij}(\sigma).$$

Let $T_{1j}$, $T_{2j}$ and $T_{3j}$ be the three terms on the last line. Using (9) and (10), we see that
(13)
$$\mathbb{E}(f(\sigma)^2) = \frac{1}{N} \sum_{j=1}^N (T_{1j} + T_{2j} + T_{3j}) m_j(\sigma) \sigma_j p_j(\sigma).$$

Now, since $\sum_i a_i(\sigma)^2 \leq 4N$ and
$$\sum_j m_j(\sigma)^2 p_j(\sigma)^2 \leq \sum_j m_j(\sigma)^2 = \|J\sigma\|^2 \leq N,$$

it follows that

$$\left|\frac{1}{N}\sum_{j=1}^{N} T_{1j} m_j(\sigma)\sigma_j p_j(\sigma)\right| = \frac{2}{N^2}\left|\sum_{i,j=1}^{N} J_{ij} a_i(\sigma) m_j(\sigma) p_j(\sigma)\right|$$

$$\leq \frac{2}{N^2}\sqrt{4N}\sqrt{N} = \frac{4}{N}.$$

Note that it would be inefficient to simply use the Cauchy–Schwarz inequality, because it would not allow us to take advantage of the large amounts of "cancellation" between positive and negative terms in the quadratic form.

Next, let us look at the $T_2$-term. We have

$$\frac{1}{N}\sum_{j=1}^{N} T_{2j} m_j(\sigma)\sigma_j p_j(\sigma) = \frac{2}{N^2}\sum_{j=1}^{N} m_j(\sigma)^2 p_j(\sigma)$$

$$\leq \frac{2}{N^2}\sum_{j=1}^{N} m_j(\sigma)^2 \leq \frac{2}{N}.$$

Finally, let us bound the $T_3$-term. Take any $i$ and let $e_i$ be the $i$th coordinate vector in $\mathbb{R}^N$. Then

$$\sum_{j=1}^{n} J_{ij}^2 = \|e_i^t J\|^2 \leq \|e_i\|^2 \|J\|^2 \leq 1.$$

Thus, if we let $J_2$ be the matrix $(J_{ij}^2)_{1\leq i,j\leq N}$, then by the well-known result that the $L^2$ operator norm of a symmetric matrix is bounded by its $L^\infty$ operator norm, we get

$$\|J_2\| \leq \max_{1\leq i\leq N}\sum_{j=1}^{N} J_{ij}^2 \leq 1.$$

Now let $h(x) := \tanh(\beta x)$. It is easy to verify that $\|h''\|_\infty \leq \beta^2$. Therefore,

$$|h(m_i(\sigma)) - h(m_i(\sigma^{(j)})) - (m_i(\sigma) - m_i(\sigma^{(j)}))h'(m_i(\sigma))|$$

$$\leq \frac{\beta^2}{2}(m_i(\sigma) - m_i(\sigma^{(j)}))^2.$$

Let $c_i(\sigma) := h'(m_i(\sigma))$, and note that $m_i(\sigma) - m_i(\sigma^{(j)}) = 2J_{ij}\sigma_j$. So the above inequality can be rewritten as

(14) $$|b_{ij}(\sigma) - 2J_{ij}\sigma_j c_i(\sigma)| \leq 2\beta^2 J_{ij}^2.$$

Finally, note that $|c_i(\sigma)| \leq \beta$. Using all this information and the bounds on the operator norms of $J$ and $J_2$, we see that for any $x,y \in \mathbb{R}^n$,

$$\left|\sum_{i,j} x_i y_j b_{ij}(\sigma)\right| \leq \left|\sum_{i,j} x_i y_j (2J_{ij}\sigma_j c_i(\sigma))\right|$$

$$+ \left| \sum_{i,j} x_i y_j (b_{ij}(\sigma) - 2 J_{ij} \sigma_j c_i(\sigma)) \right|$$

$$\leq 2 \left( \sum_i (x_i c_i(\sigma))^2 \right)^{1/2} \left( \sum_j (y_j \sigma_j)^2 \right)^{1/2} + \sum_{i,j} |x_i y_j| 2\beta^2 J_{ij}^2$$

$$\leq (2\beta + 2\beta^2) \|x\| \|y\|.$$

Again, it is clear from the definition (12) of $b_{ij}$ that $|b_{ij}(\sigma)| \leq 2\beta |J_{ij}|$. Thus,

$$\left| \sum_{i,j} x_i y_j J_{ij} b_{ij}(\sigma) \right| \leq \sum_{i,j} |x_i y_j| 2\beta J_{ij}^2 \leq 2\beta \|x\| \|y\|.$$

Applying these inequalities to the $T_3$-term in (13), we get

$$\left| \frac{1}{N} \sum_{j=1}^N T_{3j} m_j(\sigma) \sigma_j p_j(\sigma) \right| = \left| \frac{1}{N^2} \sum_{i,j=1}^N m_i(\sigma^{(j)}) b_{ij}(\sigma) m_j(\sigma) \sigma_j p_j(\sigma) \right|$$

$$= \left| \frac{1}{N^2} \sum_{i,j=1}^N (m_i(\sigma) - 2 J_{ij} \sigma_j) b_{ij}(\sigma) m_j(\sigma) \sigma_j p_j(\sigma) \right|$$

$$\leq \frac{(2\beta + 2\beta^2) + 4\beta}{N}.$$

Thus, we have computed upper bounds for all terms in (13). Combining, we have

(15) $$\mathbb{E}(f(\sigma)^2) \leq \frac{6 + 6\beta + 2\beta^2}{N}.$$

Since $S_\sigma(\beta) = f(\sigma)$, this completes the proof. □

This completes step (1) of the proof. Let us now engage in completing step (2), that is, showing that $S_\sigma(\beta) \approx 0 \Rightarrow \beta \approx \hat\beta$.

It is not difficult to verify that the function $S_\sigma$ is nonincreasing; therefore, we only have to show that $S_\sigma$ is sufficiently "nonflat" with high probability. Again, the proof is divided into two substeps: (2a) Show that if $N^{-1}H(\sigma) \geq c > 0$ for some suitable constant $c$, then $S_\sigma$ is sufficiently nonflat. (2b) Find $c > 0$ such that $N^{-1}H(\sigma) \geq c$ with high probability.

LEMMA 2.2. *Suppose $\tau \in \{-1, 1\}^N$ and $c > 0$ are such that $N^{-1}H(\tau) \geq c$. Then for any $\beta \geq 0$, we have*

$$|\tanh(2\hat\beta(\tau)/c) - \tanh(2\beta/c)| \leq \frac{8|S_\tau(\beta)|}{3c^5}.$$

PROOF. Recall the definition (5) of the local fields $m_i, i = 1, \ldots, N$. First, note that

$$N^{-1} \sum_{1 \leq i \leq N} |m_i(\tau)| \geq N^{-1} \left| \sum_{1 \leq i \leq N} m_i(\tau)\tau_i \right| = 2N^{-1}|H(\tau)| \geq 2c.$$

Also by assumption (8), we have

$$N^{-1} \sum_i m_i(\tau)^2 = N^{-1}\|J\sigma\|^2 \leq 1.$$

Now recall the Paley–Zygmund second moment inequality (see, e.g., (2.18) in [25]), which says that for any nonnegative square-integrable random variable $X$ we have

$$\mathbb{P}\{X > a\} \geq \frac{(\mathbb{E}(X) - a)^2}{\mathbb{E}(X^2)} \qquad \text{for each } 0 \leq a < \mathbb{E}(X).$$

Thus,

$$\frac{\#\{i : |m_i(\tau)| > c\}}{N} \geq \frac{(N^{-1}\sum_i |m_i(\tau)| - c)^2}{N^{-1}\sum_i m_i(\tau)^2} \geq c^2.$$

Again, by Chebyshev's bound we have

$$\frac{\#\{i : |m_i(\tau)| \geq 2/c\}}{N} \leq \frac{N^{-1}\sum_i m_i(\tau)^2}{(2/c)^2} \leq \frac{c^2}{4}.$$

Combining, we get

$$\frac{\#\{i : c < |m_i(\sigma)| < 2/c\}}{N} \geq \frac{3c^2}{4}.$$

Thus, for every $x \geq 0$,

$$|S_\tau'(x)| = \frac{1}{N} \sum_{1 \leq i \leq N} \frac{m_i(\tau)^2}{\cosh^2(xm_i(\tau))}$$

$$\geq \frac{1}{N} \sum_{i : c < |m_i(\tau)| < 2/c} \frac{m_i(\tau)^2}{\cosh^2(xm_i(\tau))}$$

$$\geq \frac{3c^2}{4} \min_{c < y < 2/c} \frac{y^2}{\cosh^2(xy)}.$$

Now, the map $y \mapsto y^2/\cosh^2(xy)$ is unimodal, and hence the minimum in the last expression is attained at either $y = c$ or $y = 2/c$. Surprisingly, the crude bound obtained by putting $y = c$ in the numerator and $y = 2/c$ in the denominator can be shown to be the best that one can do in this situation. Thus, we have

$$|S_\tau'(x)| \geq \frac{3c^4}{4\cosh^2(2x/c)}.$$

Since $S_\tau(0) = 2N^{-1}H(\tau) \geq c > 0$, $S_\tau(\infty) = N^{-1}\sum_i(m_i(\tau)\tau_i - |m_i(\tau)|) \leq 0$, and $S_\tau$ is continuous on $[0, \infty]$, therefore $S_\tau(\hat{\beta}(\tau)) = 0$ whether $\hat{\beta}(\tau) < \infty$ or not. Using this, and the monotonicity of $S_\tau$, we see that for any $\beta \geq 0$,

$$|S_\tau(\beta)| = |S_\tau(\beta) - S_\tau(\hat{\beta}(\tau))|$$
$$\geq \int_{\beta \wedge \hat{\beta}(\tau)}^{\beta \vee \hat{\beta}(\tau)} \frac{3c^4}{4\cosh^2(2x/c)}\,dx$$
$$= \frac{3c^5}{8}|\tanh(2\hat{\beta}(\tau)/c) - \tanh(2\beta/c)|.$$

This completes the proof. □

Let us now carry out step (2b)—that is, find a positive constant $c$ such that $N^{-1}H(\sigma) > c$ with high probability. The relevant result is Lemma 2.5, but we need two preliminary lemmas to overcome the problems at critical temperatures.

LEMMA 2.3. *For any nondecreasing function $f: \mathbb{R} \to \mathbb{R}$, the quantity $\mathbb{E}_\beta f(H(\sigma))$ is a nondecreasing function of $\beta$.*

PROOF. Recall the well-known inequality (e.g. [16]) that for any random variable $X$ and any two nondecreasing functions $f$ and $g$, $\text{cov}(f(X), g(X)) \geq 0$. A direct computation shows that

$$\frac{d}{d\beta}\mathbb{E}_\beta f(H(\sigma)) = \text{cov}_\beta(f(H(\sigma)), H(\sigma)).$$

This completes the proof. □

LEMMA 2.4. *For any $0 < \beta_1 \leq \beta$ and $c < \psi'(\beta_1)$, we have*

$$\mathbb{P}_\beta\{N^{-1}H(\sigma) \leq c\} \leq \frac{\psi''(\beta_1)}{N(\psi'(\beta_1) - c)^2}.$$

PROOF. The assumption that $J$ is not the zero matrix implies that $\psi(\beta)$ and $\psi'(\beta)$ are positive for all $\beta > 0$. First, assume $\beta_1 = \beta$. Standard calculations give

$$\psi'(\beta) = N^{-1}\mathbb{E}_\beta H(\sigma) \quad \text{and} \quad \psi''(\beta) = N^{-1}\text{var}_\beta H(\sigma).$$

A simple application of Chebyshev's inequality completes the proof in this case. In general, since $\beta_1 \leq \beta$, therefore Lemma 2.3 gives

$$\mathbb{P}_\beta\{N^{-1}H(\sigma) \leq c\} \leq \mathbb{P}_{\beta_1}\{N^{-1}H(\sigma) \leq c\}.$$

This completes the proof. □

We are now ready to state the step (2b) and finish the proof using the two preceding lemmas and an additional argument which eliminates the need to consider the derivatives of $\psi$.

LEMMA 2.5. *For any $\beta > 0$, we have*
$$\mathbb{P}_\beta\left\{N^{-1}H(\sigma) \leq \frac{\psi(\beta)}{4\beta}\right\} \leq \frac{8\beta^2}{N\psi(\beta)^3}.$$

PROOF. By assumption (8), we have $|H(\sigma)| = \frac{1}{2}|\sigma^t J\sigma| \leq \frac{1}{2}\|J\|\|\sigma\|^2 \leq \frac{1}{2}N$. Thus for any $\beta$,
$$0 \leq \psi'(\beta) = N^{-1}\mathbb{E}_\beta H(\sigma) \leq \tfrac{1}{2}.$$
Thus, for any $0 < \varepsilon < \beta$,
$$\psi(\beta - \varepsilon) \geq \psi(\beta) - \frac{\varepsilon}{2}.$$
Now, since $\psi'$ is an increasing function, therefore
$$(16) \qquad \psi'(\beta - \varepsilon) \geq \frac{\psi(\beta - \varepsilon) - \psi(0)}{\beta - \varepsilon} = \frac{\psi(\beta - \varepsilon)}{\beta - \varepsilon} \geq \frac{\psi(\beta) - (1/2)\varepsilon}{\beta - \varepsilon}.$$
By the mean value theorem, there exists $\beta_1 \in (\beta - \varepsilon, \beta)$ such that
$$\psi''(\beta_1) = \frac{\psi'(\beta) - \psi'(\beta - \varepsilon)}{\varepsilon}.$$
Combining this with the upper bound on $\psi'(\beta)$ and the lower bound on $\psi'(\beta - \varepsilon)$ obtained above, we get
$$\psi''(\beta_1) \leq \frac{(1/2)\beta - \psi(\beta)}{\varepsilon(\beta - \varepsilon)}.$$
Note that the numerator is always positive, since $\psi(\beta) = \int_0^\beta \psi'(u)\,du \leq \frac{1}{2}\beta$. Last of all, observe that $\psi'(\beta_1) \geq \psi'(\beta - \varepsilon)$ by the monotonicity of $\psi'$. Thus, if $\varepsilon$ and $c$ satisfy
$$(17) \qquad 0 < \varepsilon < \beta, \qquad \psi(\beta) > \frac{\varepsilon}{2} \quad \text{and} \quad c \leq \frac{1}{2}\psi'(\beta_1),$$
then Lemma 2.4 gives
$$\mathbb{P}_\beta\{N^{-1}H(\sigma) < c\} \leq \frac{4\psi''(\beta_1)}{N\psi'(\beta_1)^2}$$
$$\leq \frac{4((1/2)\beta - \psi(\beta))(\beta - \varepsilon)}{N\varepsilon(\psi(\beta) - (1/2)\varepsilon)^2}$$
$$\leq \frac{2\beta^2}{N\varepsilon(\psi(\beta) - (1/2)\varepsilon)^2}.$$

The proof is completed by putting $\varepsilon = \psi(\beta)$ and $c = \psi(\beta)/4\beta$. We only have to verify that the chosen $\varepsilon$ and $c$ satisfy (17). First, note that $\varepsilon \leq \beta/2 < \beta$, and $\psi(\beta) - \frac{1}{2}\varepsilon = \frac{1}{2}\psi(\beta) > 0$. Also, by (16),

$$c = \frac{\psi(\beta) - (1/2)\varepsilon}{2\beta} \leq \frac{\psi(\beta) - (1/2)\varepsilon}{2(\beta - \varepsilon)} \leq \frac{1}{2}\psi'(\beta - \varepsilon) \leq \frac{1}{2}\psi'(\beta_1).$$

This completes the proof. $\square$

Combining steps (2a) and (2b), we are now ready to finish the whole of step (2). The following lemma shows that $S_\sigma(\beta) \approx 0 \Rightarrow \beta \approx \hat{\beta}$.

LEMMA 2.6. *For any $\beta > 0$ and $\varepsilon > 0$, we have*

$$\mathbb{P}_\beta\{|\tanh(C_1\hat{\beta}) - \tanh(C_1\beta)| > \varepsilon\} \leq \frac{C_2}{N} + \mathbb{P}_\beta\{|S_\sigma(\beta)| > C_3\varepsilon\},$$

*where $C_1 = 8\beta/\psi(\beta)$, $C_2 = 8\beta^2/\psi(\beta)^3$ and $C_3 = 3\psi(\beta)^5/(2^{13}\beta^5)$.*

PROOF. From Lemma 2.2, we see that for any positive $\beta, \varepsilon$ and $c$,

$$\{\tau : |\tanh(2\hat{\beta}(\tau)/c) - \tanh(2\beta/c)| > \varepsilon\}$$
$$\subseteq \{\tau : N^{-1}H(\tau) \leq c\} \cup \{\tau : |S_\tau(\beta)| \geq 3c^5\varepsilon/8\}.$$

Putting $c = \psi(\beta)/(4\beta)$ and applying the bound from Lemma 2.5 completes the proof. $\square$

Finally, combining Lemmas 2.6 and 1.2, and getting rid of assumption (8) by substituting $\beta\|J\|$ for $\beta$ and $\hat{\beta}\|J\|$ for $\hat{\beta}$, we get Theorem 2.1.

**Acknowledgments.** The author is grateful to the referee and the Associate Editor for a deft understanding of the issues, and providing very useful comments.

Department of Statistics
University of California
Berkeley, California 94720-3860
USA
E-mail: sourav@stat.berkeley.edu
URL: http://www.stat.berkeley.edu/~sourav